\documentclass[12pt]{article}
\usepackage{color}
\usepackage{amscd,amsfonts,amssymb,amscd,amsmath,latexsym}

\usepackage[hypertex]{hyperref}
\usepackage[all]{xy}
\usepackage{mathrsfs}

\title{Chern classes on differential $K$-theory}
\author{Ulrich Bunke\thanks{NWF I - Mathematik,
Universit{\"a}t Regensburg,
93040 Regensburg,
GERMANY, ulrich.bunke@mathematik.uni-regensburg.de} }

\newtheorem{theorem}{Theorem}[section] 

\newtheorem{lem}[theorem]{Lemma}

\newcommand{\Rham}{{\tt Rham}}

\newcommand{\colim}{{\tt colim}}

\newcommand{\Z}{\mathbb{Z}}

\newcommand{\proof}{{\it Proof.$\:\:\:\:$}}

\newcommand{\R}{\mathbb{R}}

\newcommand{\Q}{\mathbb{Q}}

\newcommand{\bK}{{\bf K}}

\newcommand{\cE}{\mathcal{E}}

\newcommand{\cK}{\mathcal{K}}

\newcommand{\Hom}{{\tt Hom}}

\newcommand{\id}{{\tt id}}

\newcommand{\nat}{\mathbb{N}}

\def\imath{{i}}

\def\hB{\hspace*{\fill}$\Box$ \newline\noindent}

\newcommand{\ind}{{\tt index}}

\def\hB{\hspace*{\fill}$\Box$ \\[0.5cm]\noindent}

\newcommand{\bH}{\mathbf{H}}

\newcommand{\pr}{{\tt pr}}

\newcommand{\ch}{{\mathbf{ch}}}

\newcommand{\bV}{\mathbf{V}}

\newcommand{\hHZ}{{\widehat{H\Z}}}

\begin{document}\maketitle
\begin{abstract}
In this note we give a simple, model-independent  construction of  Chern classes as natural transformations from differential complex $K$-theory to differential integral cohomology. We  verify the expected behaviour of these Chern classes with respect to sums and suspension.
\end{abstract}

\tableofcontents

\section{Statements}

Complex $K$-theory and integral cohomology $H\Z$ are  generalised cohomology theories which have a unique differential\footnote{In our previous work instead of "differential cohomology" we used the term "smooth cohomology".  We were convinced by D. Freed that differential cohomology is the better name.} extensions $(\hat K,R,I,a,\int)$ and $(\hHZ,R,I,a,\int)$ with integration. 
 Moreover, these extensions are multiplicative in a unique way.  We refer to  \cite{bs2009} for a description of the axioms for differential extensions of cohomology theories and a proof of these statements.

The $i$'th Chern class is a natural transformation of \textit{set-valued} functors
$$c_i:K^0\to H\Z^{2i}$$
on the category of topological spaces.
The product 
$H\Z^{ev}:=\prod_{i\ge 0} H\Z^{2i}$ is a  functor with values in commutative graded rings.  We consider subfunctor 
$H\Z^{ev,*}_1:=1+\prod_{i\ge 1} H\Z^{2i}\subseteq \prod_{i\ge 0} H\Z^{2i}$ which takes values in the subgroup of units. 
The total Chern class
$$c:=1+c_1+c_2+\dots:K^0\longrightarrow H\Z^{ev,*}_1$$
is a natural transformation of \textit{group-valued} functors.

Let $\Omega^*_{cl}(\dots,K^*)\subseteq \Omega^*(\dots,K^*)$ denote the graded ring valued functors  on smooth manifolds of smooth differential
 forms  with coefficients in $K^*$and its subfunctor of closed forms.
We use the powers of the Bott element in $K^2$ in order to identify the functors
$$\Omega^0(\dots,K^*)\cong \Omega^{ev}(\dots) \ ,\quad  
\Omega^{-1}(\dots,K^*)\cong \Omega^{odd}(\dots)\ .$$
We therefore have  natural transformations
$$a:\Omega^{odd}\to \hat K^0\ ,\quad R:\hat K^0\to \Omega_{cl}^{ev}\ ,$$
where $a$ only preserves the additive structure, while $R$ is multiplicative.

We consider the symmetric formal power series in infinitely many variables
$$\tilde \ch:=\sum_{i\ge 1} (e^{x_i}-1)\in \Q[[x_1,x_2,\dots]]\ .$$
We write
$\ch_i$ for the homogeneous component of degree $i$.
Then there are polynomials
$$C_i\in \Q[s_1,s_2,\dots]$$ of degree $i$ (where $s_i$ has degree $i$)
such that
$$C_i(\ch_1,\dots,\ch_i)=\sigma_i$$
is the $i$th elementary symmetric function in the $x_i$.
The polynomial $C_i$ induces a natural transformation
$$C_i:\Omega^{ev}\to \Omega^{2i}$$ which maps the even form  $\omega=\omega_0+\omega_2+\omega_4+\dots$, $\omega_{2k}\in \Omega^{2k}(M)$ to
$$C_i(\omega):=C_i(\omega_2,\dots,\omega_{2i})\in \Omega^{2i}(M)\ .$$
The following theorem states that the Chern classes have unique lifts to the differential extensions which are, in addition, compatible with the group structures.

\begin{theorem}\label{th1}
\begin{enumerate}
 \item\label{aaa1} For every $i\ge 1$ there exists a unique natural transformation of set-valued functors on smooth manifolds
$$\hat c_i:\hat K^0\to \hHZ^{2i}$$ such that the following diagram commutes:
\begin{equation}\label{stern}
\xymatrix{\Omega^{ev}\ar[r]^{C_{i}}&\Omega^{2i}\\\hat K^0\ar[u]^R\ar[d]^I\ar[r]^{\hat c_i}&\hHZ^{2i}\ar[u]^R\ar[d]^I\\
K^0\ar[r]^{c_i}&H\Z^i}
\end{equation}
\item\label{aaa2} The total class
$$\hat c=1+\hat c_1+\dots:\hat K^0\to \hHZ^{ev,*}_1$$
preserves the group structure. 
\end{enumerate}
\end{theorem}

Lifts of the Chern classes have previously been constructed in \cite{berthomieu-2008}.
The goal of the present paper is to give a much simpler,  model-independent treatment.  
Further new, but not very deep, points of the present theorem are the assertions about uniqueness and the second statement. Our method of proof is different from  \cite{berthomieu-2008}. It is in fact a specialisation of
a general principle already used in \cite{bs2009} and \cite{bunke-20091} for the construction of lifts of natural transformations between cohomology functors to their differential refinements.

In the next two paragraphs we connect the differential Chern classes on differential $K$-theory  with previous constructions of differential
Chern classes in specific geometric situations.

If $\bV:=(V,h^V,\nabla^V)$ is a hermitian vector bundle with connection over a manifold $M$, then we have the Cheeger-Simons classes
$$\hat c^{CS}_i(\bV)\in \hHZ^{2i}(M)$$
constructed in \cite{MR827262}. In the model of differential $K$-theory \cite{bunke-20071} the geometric bundle is a cycle for a
differential $K$-theory class $[\bV]\in \hat K^0(M)$. We have
$$\hat c_i([\bV])=\hat c^{CS}_i(\bV)\ .$$

An even geometric family $\cE$ over $M$ (see \cite{bunke-2002} for this notion) gives rise to a Bismut superconnection
$A(\cE)$ on an infinite-dimensional Hilbert space bundle $H(\cE)$ over $M$.  This superconnection
$$A(\cE)=D(\cE)+\nabla^{H(\cE)}+\mbox{\small higher terms}$$
extends the family of Dirac operators 
  $D(\cE)$. If the kernel  of $D(\cE)$ is a vector bundle, then  it has an induced metric $h^{\ker(D(\cE))}$
and connection $\nabla^{\ker(D(\cE))}$ obtained from $\nabla^{H(\cE)}$ by projection. We thus get an induced geometric bundle
$$\bH(\cE)=(\ker(D(\cE)),h^{\ker(D(\cE))},\nabla^{\ker(D(\cE))})$$ and can define the class $\hat c^{CS}_i(\bH(\cE))\in \hHZ^{2i}(M)$. One of the original  goals of \cite{bunke-2002},
which was not quite achieved there, was to extend this construction to the general case where we do not have a kernel bundle.
Under the assumption that  $\ind(D(\cE))\in K^0(M)$ belongs to the $i$'th-step of the Atiyah-Hirzebruch filtration (i.e. vanishes
after pull-back to any $i-1$-dimensional complex) in \cite[4.1.19]{bunke-2002} we constructed a class
$\hat c_i(\cE)\in \hHZ^{2i}(M)$\footnote{Note that in \cite{bunke-2002} we index the Chern classes by their degree, where in the present note we adopt the usual convention.} such that $I(\hat c_i(\cE))=c_i(\ind(D(\cE)))$. On the other hand, the geometric family $\cE$ represents a differential $K$-theory class $[\cE,0]\in \hat K^0(M)$ in the model \cite{bunke-20071}, and we have
$I([\cE,0])=\ind(D(\cE))$.  The class $\hat c_i([\cE,0])\in \hHZ^{2i}(M)$ also satisfies
$I(\hat c_i(\cE))=c_i(\ind(D(\cE)))$ and thus
gives a second differential refinement of the $i$'th Chern class of the index of $D(\cE)$.
But in general the class
$\hat c_i(\cE)$ differs from $\hat c_i([\cE,0])$. This can already be seen on the level of curvatures. Namely, we have
$$R(\hat c_i(\cE))=R([\cE,0])_{[2i]}\ , \quad R(\hat c_i([\cE,0]))=C_i(R([\cE,0]))\ ,$$
where $\omega_{[2i]}$ denotes the degree-$2i$ component of the form $\omega$.
In a sense, the present note gives the right answer to the problem considered in \cite{bunke-2002}.

Finally we discuss odd Chern classes.
In topology, the odd Chern classes $c^{odd}_i:K^{-1}\to H\Z^{i}$ are related with the even Chern classes by suspension
$$\xymatrix{\tilde K^0(\Sigma M_+)\ar[r]^{c_{\frac{i+1}{2}}}\ar[d]^\cong&\widetilde{ H\Z}^{i+1}(\Sigma M_+)\ar[d]^\cong\\
 K^{-1}(M)\ar[r]^{c^{odd}_i}& H\Z^{i}(M)}\ .$$
In the smooth context the suspension isomorphism is replaced by the integration
$\int$ along $S^1\times M\to M$. We have the following odd counterpart of Theorem \ref{th1}.

\begin{theorem}\label{feiwfwefqfefewfqeffqewf}
For odd $i\in \nat$ there are unique natural transformations
$\hat c_i^{odd}:\hat K^{-1}\to \hHZ^{i}$ such that
$$\xymatrix{\hat  K^0(S^1\times M)\ar[r]^{\hat c_{\frac{i+1}{2}}}\ar[d]^\int& \hHZ^{i+1}(S^1\times M)\ar[d]^\int\\
\hat K^{-1}(M)\ar[r]^{\hat c^{odd}_i}& \hHZ^{i}(M)}$$ commutes.
The transformation in addition satisfies
$$I\circ \hat c^{odd}_i=c^{odd}_i\circ I\ .$$
\end{theorem}

Let $\pi:W\to B$ be a proper $K$-oriented map between manifolds. Then we have an Umkehr map
$\pi_!:K^*(W)\to K^{*-n}(B)$, where $n=\dim(W)-\dim(B)$. An integral index theorem is an assertion
about the Chern classes $c_*(\pi_!(x))$, or $c_*^{odd}(\pi_!(x))$ for $x\in K^*(W)$, e.g. an expression of these classes in terms of the classes $c_*(x)$ or $c_*^{odd}(x)$, respectively.  A prototypical example is given in \cite{madsen}.
The construction of differential lifts of Chern classes makes it possible to ask for geometric refinements of these kinds of results.
An example of such a theorem related to the Pfaffian bundle will be discussed in a forthcoming paper.

\section{Proofs}

\proof
Let $\bK_0\simeq \Z\times BU$ be a representative of the homotopy type of the classifying space of the functor $K^0$.
We choose by  \cite[Prop 2.1]{bs2009} a sequence of manifolds $(\cK_k)_{k\ge 0}$ together with maps $$x_k:\cK_k\to \bK_0\ ,\quad \kappa_k:\cK_k\to \cK_{k+1}$$ such that
\begin{enumerate}
\item $\cK_k$ is homotopy equivalent to an $i$-dimensional $CW$-complex, 
\item $\kappa_k:\cK_k\to \cK_{k+1}$ is an embedding of a closed submanifold,
\item $x_k:\cK_k\to \bK_0$ is $k$-connected,
\item  $x_{k+1}\circ \kappa_k=x_k$\ .
\end{enumerate}
Let $u\in K^0(\bK_0)$ the universal class represented by the identity map $\bK_0\to \bK_0$. By \cite[Prop. 2.6]{bs2009} we can further choose a sequence $\hat u_k\in \hat K^0(\cK_k)$ such that
$I(\hat u_k)=x_k^*u$ and 
$\kappa_k^*\hat u_{k+1}=\hat u_k$ for all $k\ge 0$.
By \cite[Lem. 3.8]{bs2009} and $2j-1<k$ we have that $H^{2j-1}(\cK_k,\R)=0$.
We consider the canonical natural transformation $\iota_\R:H\Z^*\to H\R^*$  and the de Rham map
$\Rham:\Omega^*_{cl}\to  H\R^*$.
Since $\Rham$ is multiplicative we have  $$\iota_{\R}(c_i(I(\hat u_k)))=C_i(\ch(I(\hat u_k)))=C_i(\Rham(R(\hat u_k)))=\Rham(C_i(R(\hat u_k)))\ .$$ If we choose $k\ge 2i$, then the diagram
$$\xymatrix{\hHZ^{2i}(\cK_k)\ar[r]^I\ar[d]^R&H\Z^{2i}(\cK_k)\ar[d]^{\iota_\R}\\\Omega_{cl}^{2i}(\cK_k)\ar[r]^{\Rham}&H\R^{2i}(\cK_k)}$$
is cartesian. Hence for $k\ge 2i$
there exists a unique class
$\hat z_{i,k}\in \hHZ^{2i}(\cK_k)$ such that
$$I(\hat z_{i,k})=c_i(I(\hat u_k)) \ ,\quad R(\hat z_i)=C_i(R(\hat u_k ))\ .$$ Furthermore, we have $\kappa_k \hat z_{i,k+1}=\hat z_{i,k}$. For $k<2i$ we define
$z_{i,k}:=(\kappa_k^*\circ \dots \circ \kappa_{2i-1}^* )z_{i,2i}$.

We now define the natural transformation
$\hat c_i$. We start with the observation that if $\hat c_i$ exists, then it satisfies
$$\hat c_i(\hat u_k)=\hat z_{i,k}\ .$$

Let $\hat w\in \hat K^0(M)$. 
By  \cite[Prop. 2.6]{bs2009} we have $K^0(M)\cong \colim_{k} [M,\cK_k]$, and
the underlying class $I(\hat w)\in K^0(M)$ can be written
as $I(\hat w)=f^* x_k^*u$ for some $k$ and $f:M\to \cK_k$. We choose a form $\rho\in \Omega^{odd}(M)$ such that $$\hat w=f^* \hat u_k+a(\rho)\ .$$

We consider a form
$\tilde \rho\in \Omega^{odd}([0,1]\times M)$ which restricts to $\rho$ on
$\{1\}\times M$ and to $0$ on $\{0\}\times M$.
We get a class $\tilde{\hat w}=\pr_M^*\hat w+a(\tilde \rho)\in \hat K^0([0,1]\times M)$.
Note that
$$\tilde {\hat w}_{|\{0\}\times M}=f^*\hat u_k\ ,\quad \tilde{\hat w}_{|\{1\}\times M}=\hat w\ .$$ If $\hat c_i$ exists, then we must have by naturality and the homotopy formula \cite[(1)]{bs2009}
$$\hat c_i(\tilde {\hat w}_{|\{0\}\times M})=f^*\hat z_{i,k}\ ,\quad \hat c_i(
 \tilde{\hat w}_{|\{1\}\times M}) -\hat c_i(\tilde {\hat w}_{|\{0\}\times M})=a(\int_{[0,1]\times M/M}  R(\hat c_i(\tilde{\hat w})))\ .$$
Furthermore, by the commutativity of the upper square in (\ref{stern}) we  must require
$$R(\hat c_i(\tilde {\hat w}))=C_i(R(\tilde {\hat w})) \ .$$
 Therefore we are forced to define
\begin{equation}\label{eq1}
\hat c_i(\hat w):=f^*\hat z_{i,k}+a(\int_{[0,1]\times M/M}   C_i(R(\tilde {\hat w})) )
\end{equation}
We see that if $\hat c_i$ exists, then it is automatically unique.
\begin{lem}
The definition of  $\hat c_i(\hat w)$ by (\ref{eq1}) is independent of the choices of $\tilde \rho$, $\rho$ and $f:M\to \cK_k$.
\end{lem}
\proof
Let us start with a second choice $\tilde \rho^\prime$ and write $\tilde{\hat w}^\prime:=\pr_M^*\hat w+a(\tilde \rho^\prime)$. Then we can connect $\tilde \rho$ with $\tilde \rho^\prime$ by a family of such forms, e.g. the linear path. This path can be considered as a form $\bar \rho$ on $[0,1]\times [0,1]\times M$.  By construction
$\bar \rho_{|[0,1]\times \{j\}\times M}$ is constant and has no component in the
direction of the first variable for $j=0,1$ .
This implies that
\begin{equation}\label{eq2}
R(\tilde {\hat w}^\prime)_{|[0,1]\times \{j\}\times M}=0\ .
\end{equation}
We set $\bar{\hat w}:=\pr_M^*\hat w+a(\bar \rho)\in \hat K^0([0,1]\times [0,1]\times M)$.
By Stokes theorem we have 
$$d\int_{[0,1]\times [0,1]\times M/M} C_i(R((\bar{\hat w})))=
\int_{[0,1]\times M/M}   C_i(R(\tilde {\hat w}^\prime)) - \int_{[0,1]\times M/M}   C_i(R(\tilde {\hat w}))$$
(these are the contributions of the faces $\{j\}\times [0,1]\times M$)
since the integral over the other two faces $[0,1]\times \{j\}\times M$ vanishes by (\ref{eq2}).
Since $a$ annihilates exact forms this implies that
$$a(\int_{[0,1]\times M/M}   C_i(R(\tilde {\hat w})) )=a(\int_{[0,1]\times M/M}   C_i(R(\tilde {\hat w}^\prime)) )\ .$$

Assume now that we have chosen a different $\rho^\prime$. 
Then $a(\rho^\prime-\rho)=0$ so that by the exactness axiom \cite[(2)]{bs2009} there exists a class
$\hat v\in \hat K^1(M)$ with $R(\hat v)=\rho^\prime-\rho$.
Let $\hat e\in \hat K^1(S^1)$ be a lift of the generator of $K^1(S^1)\cong \Z$ with $R(\hat e)=dt$.
We consider the form $\tilde \sigma\in \Omega^{odd}([0,1]\times M)$ with no $dt$-component given by 
$$\tilde \sigma_{|\{t\}\times M}:= \int_{[0,t]\times M/M}R(\hat e\times \hat v)\ ,$$
where we identify $S^1\cong \R/\Z$ and view the interval $[0,t]$ as a subset of $S^1$.
Then $$\tilde \sigma_{|\{0\}\times M}=0\ ,\quad \tilde \sigma_{|\{1\}\times M}=\rho^\prime-\rho¸\ ,\quad d\tilde \sigma=dt\wedge \pr_M^* R(\hat v)=R(\hat e\times \hat v)\ .$$
We now consider
$$\tilde{\hat v}:=\pr_M^*\hat w+\pr_M^*a(\rho)+a(\tilde \sigma)\in \hat K^0([0,1]\times M)$$
and calculate modulo the image of $d$
\begin{eqnarray*}
\int_{[0,1]\times M/M} C_i(R( \tilde{\hat v}))&\equiv&\int_{S^1\times M/M} C_i(R(\pr_M^* (\hat w))+\pr_M^*d \rho+R(\hat e\times \hat v))\\ 
&\equiv& \int_{S^1\times M/M} C_i(R(\pr_M^* (\hat w)) +R(\hat e\times \hat v)) \\
&\equiv& \int_{S^1\times M/M} C_i(R(\pr_M^* (\hat w) +\hat e\times \hat v))\ .
\end{eqnarray*}
It follows that
\begin{eqnarray*}
\Rham (
\int_{[0,1]\times M/M} C_i(R( \tilde{\hat v})))&=&\Rham ( \int_{S^1\times M/M} C_i(R(\pr_M^* (\hat w) +\hat e\times \hat v))) \\&=&\int_{S^1\times M/M}\Rham ( C_i(R(\pr_M^* (\hat w) +\hat e\times \hat v)))\\&=&\int_{S^1\times M/M}\iota_\R (c_i (I(\pr_M^* (\hat w) +\hat e\times \hat v)))\ .
\end{eqnarray*}
 In other words,
$\Rham (\int_{[0,1]\times M/M} C_i(R( \tilde{\hat v})))$
is  an integral class,
and this
implies
$$a(\int_{[0,1]\times M/M} C_i(R( \tilde{\hat v})))=0$$ by \cite[(2)]{bs2009}.

If $\tilde \rho$ was the path connecting $\rho$ with $0$, then we construct the path
$\tilde \rho^\prime$ from $\rho^\prime$ to $0$  by concatenating $\tilde \rho$ with $\tilde \sigma$
(in order to concatenate smoothly we can change $\tilde \rho$).
Then  get $\tilde {\hat w}^\prime:=\pr_M^*\hat w+a(\tilde \rho^\prime)\in \hat K^0([0,1]\times M)$ and 
\begin{eqnarray*}
a(\int_{[0,1]\times M/M}  C_i(R( \tilde{\hat w}^\prime)))&=&a(\int_{[0,1]\times M/M} C_i(R( \tilde{\hat w})))\\&&+a(\int_{[0,1]\times M/M} C_i(R( \tilde{\hat v})))\\&=&
a(\int_{[0,1]\times M/M} C_i(R( \tilde{\hat w})))
\end{eqnarray*}
This finishes the verification that our construction of $c_i$ is independent of the choice of $\rho$.
 
Finally we verify that $\hat c_i(\hat w)$ is independent of the choice of $f:M\to \cK_k$.
If we replace $k$ by $k+1$ and $f$ by $\kappa_k\circ f$, then we obviously get the same result. For two choices $f:M\to \cK_k$ and $f^\prime:M\to \cK_{k^\prime}$ there exists $k^{\prime\prime}\ge \max\{k,k^\prime\}$  such that
$\kappa_k^{k^{\prime\prime}}\circ f$ and $\kappa_{k^\prime}^{k^{\prime\prime}}\circ f^\prime$
are homotopic. Here $\kappa_i^j:\cK_i\to \cK_j$ denotes for $j>i$ the composition
$\kappa_i^j:=\kappa_{j-1}\circ\dots\circ \kappa_i$.  Therefore it remains to show that
a choice
$f^\prime:M\to \cK_k$ homotopic to $f:M\to \cK_k$ gives the same result for $\hat c_i(\hat w)$.
Let $H:[0,1]\times M\to \cK_k$ be a homotopy from $f$ to $f^\prime$.
Then we use $H$ in the construction of $\hat c_i(\pr_M^*\hat w)\in \hHZ^{2i}([0,1]\times M)$. If we let $\hat c_{i}^\prime(\hat w)$ denote the result of the construction based on the choice of  $f^\prime$ we have by the homotopy formula
$$\hat c_{i}^\prime(\hat w)-\hat c_{i}^\prime(\hat w)=a(\int R(\hat c_i(\pr_M^*\hat w) ))=
a(\int \pr_M^* C_i(\hat w))=0\ .$$
\hB 

\begin{lem}
The construction of $\hat c_i$ defines a natural transformation $\hat c_i:\hat K\to \hHZ^{2i}$ of set-valued  functors on smooth manifolds.
\end{lem}
\proof
Let $g:N\to M$ be a smooth map between manifolds. Let $\hat w\in \hat K^0(M)$ and assume that we have constructed
$\hat c_i(\hat w)$ using the choices of $f:M\to \cK_k$, $\rho\in \Omega^{odd}(M)$ and $\tilde \rho\in \Omega^{odd}([0,1]\times M)$.
Then we construct
$\hat c_i(g^*\hat w)$ using the choices $f\circ g:N\to \cK_k$ and
$g^*\rho\in \Omega^{odd}(N)$, $(\id\times g)^*\tilde \rho\in \Omega^{odd}([0,1]\times N)$. With these choices we have $(\id\times g)^*\tilde{\hat w}=\widetilde{g^*\hat w}\in \hat K^0([0,1]\times N)$ and 
\begin{eqnarray*}
g^*\hat c_i(\hat w)&=&g^*f^*\hat z_{i,k}+g^*a(\int_{[0,1]\times M/M} C_i(R(\tilde{\hat w})))\\&=&
(f\circ g)^*\hat z_{i,k}+ a(\int_{[0,1]\times M/M} C_i(R((\id\times g)^*\tilde{\hat w})))\\
&=&(f\circ g)^*\hat z_{i,k}+ a(\int_{[0,1]\times M/M} C_i(R(\widetilde{g^*\hat w} )))\\
&=&\hat c_i(g^*\hat w)\ .
\end{eqnarray*}
\hB
This finishes the proof of Assertion \ref{aaa1} of Theorem \ref{th1}.

In order to show the second Assertion \ref{aaa2} 
we consider the natural transformation
$$\hat B:\hat K^0\times \hat K^0\to \hHZ^{ev}$$ given by
$$\hat B(\hat w,\hat v):=\hat c(\hat w)\cup \hat c(\hat v)-\hat c(\hat w+\hat v)\in \hHZ^{ev}(M)\ , \quad \hat w,\hat v\in \hat K^0(M)\ .$$
If we apply $I$ we get
 \begin{eqnarray*}
I(\hat B(\hat w,\hat v))&=&I(\hat c(\hat w)\cup \hat c(\hat v))-I(\hat c(\hat w+\hat v))\\
&=&I(\hat c(\hat w))\cup I(\hat c(\hat v))-I(\hat c(\hat w+\hat v))\\
&=&c(I(\hat w))\cup c(I(\hat v))-c(I(\hat w)+I(\hat v))\\
&=&0\ .
\end{eqnarray*}
Let $C=1+C_1+C_2+\dots\in \Q[[s_0,s_1,\dots]]$.
Then we have the identity
$$C(s_0+s_0^\prime,s_1+s_1^\prime,\dots)=C(s_0,s_1,\dots)C(s_0^\prime,s_1^\prime,\dots)\ .$$
Indeed, if 
$$\tilde \ch=\sum_{i\ge 1} (e^{x_i}-1)\ ,$$
then
$$C(\ch_1,\dots)=\prod_{i\ge 1} (1+x_i)\ .$$
If we introduce another set of variables $x_i^\prime$ and set $\tilde \ch^\prime=\sum_{i\ge 1} (e^{x^\prime_i}-1)$,
then
\begin{eqnarray*}
C(\ch_1+\ch_1^\prime,\ch_2+\ch_2^\prime,\dots)&=&\prod_{i\ge 1} (1+x_i)(1+x_i^\prime)\\&=&
C(\ch_1,\ch_2,\dots)C(\ch_1^\prime,\ch_2^\prime,\dots)\ . 
\end{eqnarray*}
We now calculate
\begin{eqnarray*}
R(\hat B(\hat w,\hat v))&=&
R(\hat c(\hat w)\cup \hat c(\hat v))-R(\hat c(\hat w+\hat v))\\
&=&
R(\hat c(\hat w))\cup R(\hat c(\hat v))-R(\hat c(\hat w+\hat v))\\
&=&C(R(\hat w))\wedge C(R(\hat v))-C(R(\hat w)+R(\hat v))\\
&=&0\ .
\end{eqnarray*}
It follows that
$\hat B$ factorises over the subfunctor
$$H\R^{odd}/H\Z^{odd}\subset H\R/\Z^{odd}\subset \hHZ^{ev}\ ,$$
where the inclusion is induced by $a$.
Let $\rho\in \Omega^{odd}(M)$ and consider
$\tilde \rho:=t\pr_M^*\rho\in \Omega^{odd}([0,1]\times M)$.
Then we have
$$\hat B(\hat w+a(\rho),\hat v)-\hat B(\hat w ,\hat v)=\hat B(\pr_M^*\hat w+a(\tilde \rho),\hat v)_{|\{1\}\times M}-\hat B(\pr_M^*\hat w+a(\tilde \rho),\hat v)_{|\{0\}\times M}\ .$$
Since $\hat B$ takes values in the homotopy invariant subfunctor 
$H\R^{odd}/H\Z^{odd}$ we conclude that 
$\hat B(\hat w+a(\rho),\hat v)=\hat B(\hat w ,\hat v)$. In a similar manner we see that
$\hat B(\hat w,\hat v+a(\rho))=\hat B(\hat w ,\hat v)$.
Hence $\hat B$ has a factorisation over a natural transformation
$$K^0\times K^0\to H\R^{odd}/H\Z^{odd}\subset H\R/\Z^{odd}\ .$$
Such a natural transformation between homotopy invariant functors on  manifolds
must be represented by a map of classifying spaces
$$\bK_0\times \bK_0\to K(\R/\Z,odd)\ ,$$
where
$K(\R/\Z,odd):=\bigvee_{i\ge 0}K(\R/\Z,2i+1)$ is a wedge of Eilenberg-MacLane spaces,
 i.e. by  a class in $B\in H^{odd}(\bK_0\times \bK_0;\R/\Z)$.
Since $\bK_0$ and therefore $\bK_0\times \bK_0$ are even spaces
we know that
$H_{odd}(\bK_0\times \bK_0;\Z)=0$.  
It follows by the universal coefficient formula that
$H^{odd}(\bK_0\times \bK_0;\R/\Z)\cong \Hom(H_{odd}(\bK_0\times \bK_0;\Z),\R/\Z)=0$.
We see that $B=0$ and therefore $\hat B=0$.
This finishes the proof of Assertion \ref{aaa2} of Theorem \ref{th1}. \hB

We now show Theorem \ref{feiwfwefqfefewfqeffqewf}.
We let $\hat e\in K^{1}(S^1)$ be, as above,  the unique element with $R(\hat e)=dt$, $I(\hat e)=e\in K^1(S^1)$ the canonical generator, and $\hat e_{|*}=0$ for a basepoint $*\in S^1$.
Then we define for odd $i\in \nat$ and $\hat x\in \hat K^{-1}(M)$
$$\hat c_{i}^{odd}(\hat x):=\int \hat c_{\frac{i+1}{2}}(\hat e\times \hat x)\ .$$
Note that
$$I(\int\hat c_{\frac{i+1}{2}}(\hat e\times \hat x))=\int c_{\frac{i+1}{2}}(e\times I(\hat x))\ .$$
We have a natural inclusion
$\widetilde{H\Z}^{*}(\Sigma M_+)\subset H\Z^*(S^1\times M)$ as the subspace of classes whose restriction to $\{*\}\times M$ vanishes. Since $e_{|*}=0$ we see that $e\times I(\hat x)$ belongs to this subspace. The restriction of $\int$ to this subspace coincides with the suspension isomorphism 
$\widetilde{H\Z}^{*+1}(\Sigma M_+)\stackrel{\sim}{\to} H\Z^{*}(M)$, $\int(e\times x)=x$ with inverse $x\mapsto e\times x$. Therefore
$$\int c_{\frac{i+1}{2}}(e\times I(\hat x))=c_i^{odd}(I(\hat x))\ .$$
In this way we get a natural transformation which has the required property.

 Since $\int:\hat K^0(S^1\times M)\to \hat K^{-1}(M)$ is surjective it is clear that $\hat c_i^{odd}$ is unique.
\hB


\end{document}